\documentclass[12pt]{amsart}

\usepackage{amsmath, amssymb, amsthm, longtable, booktabs, url}
\usepackage{graphicx}
\usepackage{hyperref}

\newcommand{\KhSpace}{\mathcal{X}_\mathit{Kh}}
\newcommand{\St}{\mathit{St}}

\newcommand{\C}{\mathbb{C}}
\newcommand{\F}{\mathbb{F}}
\newcommand{\N}{\mathbb{N}}
\newcommand{\Z}{\mathbb{Z}}

\newcommand{\iso}{\cong}

\newcommand{\Kh}{\mathit{Kh}}

\DeclareMathOperator{\im}{im}
\DeclareMathOperator{\Sq}{Sq}

\newtheorem{theorem}{Theorem}

\newtheorem{remark}[theorem]{Remark}

\title[Computations of the Lipshitz-Sarkar Steenrod
  Square]{Computations of the Lipshitz-Sarkar Steenrod Square on
  Khovanov Homology}

\author[Cotton Seed]{Cotton Seed}
\address{Department of Mathematics \\ Princeton University}
\email{cseed@math.princeton.edu}

\begin{document}
\begin{abstract} 
  Lipshitz and Sarkar recently introduced a space-level refinement of
  Khovanov homology.  This refinement induces a Steenrod square
  operation $\Sq^2$ on Khovanov homology which they describe
  explicitly.  This paper presents some computations of $\Sq^2$.  In
  particular, we give examples of links with identical integral
  Khovanov homology but with distinct Khovanov homotopy types.
\end{abstract}

\maketitle

Let $L$ be a link in $S^3$.  Khovanov homology associates to $L$ a
bigraded abelian group $\Kh^{i,j}(L; \Z)$ \cite{khovanov}
\cite{bncat}.  Recently, Lipshitz and Sarkar defined a space-level
refinement of Khovanov homology \cite{lskhhtt}.  They construct a
family of stable spaces $\KhSpace^j(L)$ whose homotopy type is an
invariant of $L$ and whose reduced singular homology is Khovanov
homology:
\[ \widetilde{H}^i(\KhSpace^j(L); \F_2) = \Kh^{i,j}(L; \F_2). \]
The spaces $\KhSpace^j(L)$ carry stable cohomology operations.  In the
sequel \cite{lssteenrodsq}, they give an explicitly computable
combinatorial definition of the second Steenrod square operation
\[ \Sq^2 : \Kh^{i,j}(L; \F_2)\to \Kh^{i+2,j}(L; \F_2). \]
They implemented this operation in Sage and computed the operation for
all prime links with 11 or fewer crossings.  Despite these
computations, they did not find an example where $\KhSpace$
distinguishes links with identical Khovanov homology.

There is also a simpler Steenrod square
\[ \Sq^1 : \Kh^{i,j}(L; \F_2)\to \Kh^{i+1,j}(L; \F_2) \]
given by the composition of coefficient reduction $\Z\to \F_2$ with
the Bockstein homomorphism.  It is determined by $\Kh(L; \Z)$.

The author has implemented support for computing $\Sq^1$ and $\Sq^2$
in {\tt knotkit} \cite{knotkit}, a C++ software package for performing
knot homology computations, following the algorithm outlined in
\cite[Appendix A]{lssteenrodsq}.  Using this implementation, we have
computed $Sq^2$ for all knots and hyperbolic links with 14 or fewer
crossings and a large collection of 15 crossing knots.  For link data,
we use the HTW knot tables \cite{htwknots} and the Thistlethwaite link
(MT) tables and from SnapPy \cite{snappy}.  This paper presents the
results of those computations.

Lipshitz and Sarkar define a map $\St = \St(L) : \Z^2\to \N^4$ as
follows \cite[Definition 4.3]{lssteenrodsq}.  Fix $i, j\in \Z$.  For
$k \in \{i, i + 1\}$, let $\Sq^1_{(k)} = \Sq^1 : \Kh^{k,j}(L; \F_2)\to
\Kh^{k+1,j}(L; \Z)$.  Let $r_1$ be the rank of $\Sq^2 : \Kh^{i,j}(L;
\F_2)\to \Kh^{i+2,j}$, let $r_2$ be the rank of $\Sq^2|_{ker
  \Sq^1_{(i)}}$, let $r_3$ be the dimension of $\im \Sq^1_{(i+1)} \cap
\im \Sq^2$ and let $r_4$ be the dimension of $\im \Sq^1_{(i+1)} \cap
\im(\Sq^2|_{\ker \Sq^1_{(i)}})$.  Then
\[ \St(i, j) = (r_2 - r_4, r_1 - r_2 - r_3 + r_4, r_4, r_3 - r_4). \]
The map $\St(L)$ is manifestly an invariant of $\KhSpace(L)$.

First, we give a positive answer to \cite[Question 5.1]{lssteenrodsq}:
\begin{theorem}
  The space-level invariant $\KhSpace(L)$ is a strictly stronger than
  Khovanov homology.  There exist pairs of links $L$ and $L'$ such
  that $\Kh(L; \Z)\iso \Kh(L'; \Z)$ but $\KhSpace^j(L)$ and
  $\KhSpace^j(L)$ are not stably homotopy equivalent for some $j$.  A
  list of groups of such links is given in Table~\ref{table}.
\end{theorem}
All the links in Table~\ref{table} have homological width 3 over
$\F_2$.  We only list groups of links with identical integral Khovanov
homology that are distinguished by $\St(L)$.  All groups are pairs
except the lone triple $\text{K}14n5017$, $\text{K}14n11311$ and
$\text{K}14n11629$.

Khovanov homology with $\F_2$ coefficients is Conway mutation
invariant \cite{bloomodd} \cite{wehrlimutinv}.  Integral Khovanov
homology is known not to be mutation invariant for mutations that mix
components \cite{wehrlimut}, but the question of mutation invariance
for component-preserving mutations is still unaswered.  The question
of mutation invariance of $\KhSpace$ is particularly interesting.  We
have computed $\Sq^2$ for all mutant knot groups with 14 or fewer
crossings and 53,749 of the 66,108 15 crossing knots with mutant
pairs.  We used the mutant knot tables compiled by Stoimeno
\cite{stoimenow}.  The author is currently unaware of tables for
mutant links.
\begin{remark}
  For every mutant pair $L$ and $L'$ in our computation, we find
  $St(L) = St(L')$.  In particular, we find no counterexample to
  mutation invariance of $\KhSpace$.
\end{remark}

Finally, for links with suitably simple Khovanov homology, Lipshitz
and Sarkar show the homotopy type of $\KhSpace^j(L)$ is determined by
$\St(i, j) = (a, b, c, d)$ for some fixed $i$.  In this case,
$\KhSpace^j(L)$ is homotopy equivalent to a wedge sum of known spaces,
and $a$ gives the number of $\Sigma^m \C P^2$ summands in the wedge
sum decomposition.  In their computations of links with 11 or fewer
crossings, they notice that no $\C P^2$ summands are present, see
\cite[Question 5.2]{lssteenrodsq}.  This pattern persists in our
computations as well, and we make the following remark.
\begin{remark}
  For all links $L$ in our computation and for all $(i, j)$, we find
  $a = 0$ where $\St(L)(i, j) = (a, b, c, d)$.
\end{remark}

\appendix
\section{Computations}

\tiny 
\begin{center}
\begin{longtable}{p{0.13\textwidth}@{}p{0.85\textwidth}}
\caption{Groups of links with identical integral Khovanov homology
  that are distinguished by $St(L)$.}\label{table}\\*
\toprule
$L$ & $\St(L)$ \\
\midrule
$\text{K}12n41$ & $(-4, -5) \mapsto (0, 1, 0, 0)$, $(-3, -3) \mapsto (0, 0, 1, 0)$, $(-2, -1) \mapsto (0, 1, 0, 0)$, $(-1, 1) \mapsto (0, 1, 0, 0)$, $(0, 3) \mapsto (0, 0, 1, 0)$ \\
$\text{K}13n2067$ & $(-4, -5) \mapsto (0, 1, 0, 0)$, $(-3, -3) \mapsto (0, 0, 1, 0)$, $(-1, 1) \mapsto (0, 1, 0, 0)$, $(0, 3) \mapsto (0, 0, 1, 0)$ \\
\midrule
$\text{K}11n70$ & $(-3, -7) \mapsto (0, 1, 0, 0)$, $(-2, -5) \mapsto (0, 0, 1, 0)$, $(0, -1) \mapsto (0, 0, 1, 0)$ \\
$\text{K}13n2566$ & $(-3, -7) \mapsto (0, 0, 0, 1)$, $(0, -1) \mapsto (0, 0, 1, 0)$ \\
\midrule
$\text{K}13n1529$ & $(2, 9) \mapsto (0, 1, 0, 1)$, $(3, 11) \mapsto (0, 1, 1, 0)$, $(4, 13) \mapsto (0, 1, 1, 0)$, $(5, 15) \mapsto (0, 1, 1, 0)$, $(6, 17) \mapsto (0, 1, 1, 0)$, $(7, 19) \mapsto (0, 0, 1, 0)$ \\
$\text{K}14n2037$ & $(1, 7) \mapsto (0, 1, 0, 0)$, $(2, 9) \mapsto (0, 1, 1, 0)$, $(3, 11) \mapsto (0, 1, 1, 0)$, $(4, 13) \mapsto (0, 1, 1, 0)$, $(5, 15) \mapsto (0, 1, 1, 0)$, $(6, 17) \mapsto (0, 1, 1, 0)$, $(7, 19) \mapsto (0, 0, 1, 0)$ \\
\bottomrule
\end{longtable}
\end{center}

\end{document}